\documentclass[11pt]{amsart}
\usepackage{amssymb,mathrsfs}
\title{Invariant Metrics and Laplacians on Siegel-Jacobi Space}
\author{ Jae-Hyun Yang}
\address{Department of Mathematics, Inha University,
Incheon 402-751, Korea}
\email{jhyang@inha.ac.kr }

\begin{document}

\thanks{\noindent{2000 Mathematics Subject Classification:} Primary
32F45, 32M10.}
\thanks{This work was supported by Inha University Research Grant
}

\keywords{invariant metrics, Siegel-Jacobi space, spectral theory}

\maketitle
\begin{abstract}In this paper, we compute Riemannian metrics on the
Siegel-Jacobi space which are invariant under the natural action
of the Jacobi group explicitly and also provide the Laplacians of
these invariant metrics. These are expressed in terms of the trace
form.
\end{abstract}

\newtheorem{theorem}{Theorem}[section]
\newtheorem{lemma}{Lemma}[section]
\newtheorem{proposition}{Proposition}[section]
\newtheorem{remark}{Remark}[section]
\newtheorem{definition}{Definition}[section]

\renewcommand{\theequation}{\thesection.\arabic{equation}}
\renewcommand{\thetheorem}{\thesection.\arabic{theorem}}
\renewcommand{\thelemma}{\thesection.\arabic{lemma}}
\newcommand{\bbr}{\mathbb R}
\newcommand{\bbs}{\mathbb S}
\newcommand{\bn}{\bf n}
\def\charf {\mbox{{\text 1}\kern-.24em {\text l}}}

\newcommand\BD{\mathbb D}
\newcommand\BH{\mathbb H}
\newcommand\BR{\mathbb R}
\newcommand\BC{\mathbb C}
\newcommand\lrt{\longrightarrow}
\newcommand\lmt{\longmapsto}
\newcommand\CX{{\Cal X}}
\newcommand\td{\bigtriangledown}
\newcommand\pdx{ {{\partial}\over{\partial x}} }
\newcommand\pdy{ {{\partial}\over{\partial y}} }
\newcommand\pdu{ {{\partial}\over{\partial u}} }
\newcommand\pdv{ {{\partial}\over{\partial v}} }
\newcommand\PZ{ {{\partial}\over {\partial Z}} }
\newcommand\PW{ {{\partial}\over {\partial W}} }
\newcommand\PZB{ {{\partial}\over {\partial{\overline Z}}} }
\newcommand\PWB{ {{\partial}\over {\partial{\overline W}}} }
\newcommand\PX{ {{\partial}\over{\partial X}} }
\newcommand\PY{ {{\partial}\over {\partial Y}} }
\newcommand\PU{ {{\partial}\over{\partial U}} }
\newcommand\PV{ {{\partial}\over{\partial V}} }
\renewcommand\th{\theta}
\renewcommand\l{\lambda}
\renewcommand\k{\kappa}
\newcommand\G{\Gamma}
\newcommand\s{\sigma}
\newcommand\g{\gamma}
\newcommand\PWE{ \frac{\partial}{\partial \overline W}}

\begin{section}{{\bf Introduction}}
\setcounter{equation}{0}
\end{section}

For a given fixed positive integer $n$, we let
$${\mathbb H}_n=\,\{\,Z\in \BC^{(n,n)}\,|\ Z=\,^tZ,\ \ \ \text{Im}\,Z>0\,\}$$
be the Siegel upper half plane of degree $n$ and let
$$Sp(n,\BR)=\{ M\in \BR^{(2n,2n)}\ \vert \ ^t\!MJ_nM= J_n\ \}$$
be the symplectic group of degree $n$, where
$$J_n=\begin{pmatrix} 0&E_n \\
                   -E_n&0 \\ \end{pmatrix}.$$
We see that $Sp(n,\BR)$ acts on $\BH_n$ transitively by
\begin{equation}M\cdot Z=(AZ+B)(CZ+D)^{-1}, \end{equation}
where $M=\begin{pmatrix} A&B\\ C&D\end{pmatrix}\in Sp(n,\BR)$ and
$Z\in \BH_n.$

For two positive integers $n$ and $m$, we consider the Heisenberg
group
$$H_{\BR}^{(n,m)}=\{\,(\l,\mu;\k)\,|\ \l,\mu\in \BR^{(m,n)},\ \k\in \BR^{(m,m)},\ \
\k+\mu\,^t\l\ \text{symmetric}\ \}$$ endowed with the following
multiplication law
$$(\l,\mu;\k)\circ (\l',\mu';\k')=(\l+\l',\mu+\mu';\k+\k'+\l\,^t\mu'-
\mu\,^t\l').$$ We define the semidirect product of $Sp(n,\BR)$ and
$H_{\BR}^{(n,m)}$
$$G^J:=Sp(n,\BR)\ltimes H_{\BR}^{(n,m)}$$
endowed with the following multiplication law
$$\Big(M,(\lambda,\mu;\kappa)\Big)\cdot\Big(M',(\lambda',\mu';\kappa')\Big) =\,
\Big(MM',(\tilde{\lambda}+\lambda',\tilde{\mu}+ \mu';
\kappa+\kappa'+\tilde{\lambda}\,^t\!\mu'
-\tilde{\mu}\,^t\!\lambda')\Big)$$ with $M,M'\in Sp(n,\BR),
(\lambda,\mu;\kappa),\,(\lambda',\mu';\kappa') \in
H_{\BR}^{(n,m)}$ and
$(\tilde{\lambda},\tilde{\mu})=(\lambda,\mu)M'$. We call this
group $G^J$ the {\it Jacobi group} of degree $n$ and index $m$. We
have the {\it natural action} of $G^J$ on $\BH_n\times
\BC^{(m,n)}$ defined by
\begin{equation}\Big(M,(\lambda,\mu;\kappa)\Big)\cdot (Z,W)=\Big(M\cdot Z,(W+\lambda Z+\mu)
(CZ+D)^{-1}\Big), \end{equation}
where $M=\begin{pmatrix} A&B\\
C&D\end{pmatrix} \in Sp(n,\BR),\ (\lambda,\mu; \kappa)\in
H_{\BR}^{(n,m)}$ and $(Z,W)\in \BH_n\times \BC^{(m,n)}.$ The
homogeneous space $\BH_n\times \BC^{(m,n)}$ is called the {\it
Siegel}-{\it Jacobi  space} of degree $n$ and index $m.$ We refer
to [2-3], [6-7], [11], [14-21] for more details on materials
related to the Siegel-Jacobi space.

\vskip 0.31cm
\newcommand\bz{d{\overline Z}}
\newcommand\bw{d{\overline W}}
\newcommand\SJ{{\mathbb H}_n\times {\mathbb C}^{(m,n)}}
\newcommand\tr{\triangledown}
\newcommand\Hnm{{\mathbb H}_{n,m}}
\newcommand\Hn{{\mathbb H}_n}
\newcommand\Cmn{{\mathbb C}^{(m,n)}}

For brevity, we write $\Hnm:=\SJ.$ For a coordinate $(Z,W)\in\Hnm$
with $Z=(z_{\mu\nu})\in {\mathbb H}_n$ and $W=(w_{kl})\in \Cmn,$
we put
\begin{eqnarray*}
Z\,=&\,X\,+\,iY,\quad\ \ X\,=\,(x_{\mu\nu}),\quad\ \
Y\,=\,(y_{\mu\nu})
\ \ \text{real},\\
W\,=&U\,+\,iV,\quad\ \ U\,=\,(u_{kl}),\quad\ \ V\,=\,(v_{kl})\ \
\text{real},\\
dZ\,=&\,(dz_{\mu\nu}),\quad\ \ d{\overline Z} \,=\,(d{\overline
z}_{\mu\nu}),\quad\ \
dY\,=\,(dy_{\mu\nu}),\\
dW\,=&\,(dw_{kl}),\quad\ \ d{\overline W}\,=\,(d{\overline
w}_{kl}),\quad\ \ dV\,=\,(dv_{kl}),
\end{eqnarray*}

\begin{eqnarray*}
\PZ\,=\,&\left(\, { {1+\delta_{\mu\nu}} \over 2}\, {
{\partial}\over {\partial z_{\mu\nu}} } \,\right),\quad
\PZB\,=\,\left(\, { {1+\delta_{\mu\nu}}\over 2} \, {
{\partial}\over {\partial {\overline z}_{\mu\nu} }  }
\,\right),\\
\PX\,=\,&\left(\, { {1+\delta_{\mu\nu}}\over 2}\, {
{\partial}\over {\partial x_{\mu\nu} } } \,\right),\quad
\PY\,=\,\left(\, { {1+\delta_{\mu\nu}}\over 2}\, { {\partial}\over
{\partial y_{\mu\nu} }  } \,\right),
\end{eqnarray*}

$$\PW=\begin{pmatrix} {\partial}\over{\partial w_{11}} & \hdots &
 {\partial}\over{\partial w_{m1}} \\
\vdots&\ddots&\vdots\\
 {\partial}\over{\partial w_{1n}} &\hdots & {\partial}\over
{\partial w_{mn}} \end{pmatrix},\quad \PWB=\begin{pmatrix}
{\partial}\over{\partial {\overline w}_{11} }   &
\hdots&{ {\partial}\over{\partial {\overline w}_{m1} }  }\\
\vdots&\ddots&\vdots\\
{ {\partial}\over{\partial{\overline w}_{1n} }  }&\hdots &
 {\partial}\over{\partial{\overline w}_{mn} }  \end{pmatrix},$$
$$\PU=
\begin{pmatrix} {\partial}\over{\partial u_{11}} &\hdots&
{\partial}\over{\partial u_{m1}}  \\
\vdots &\ddots &\vdots\\
 {\partial}\over{\partial u_{1n}}  &\hdots &
{\partial}\over{\partial u_{mn}} \end{pmatrix},\quad \PV=
\begin{pmatrix} { {\partial}\over{\partial v_{11}} } &\hdots &
{ {\partial}\over{\partial v_{m1}} }\\
\vdots &\ddots &\hdots\\
{ {\partial}\over{\partial v_{1n}} } &\hdots &
 {\partial}\over{\partial v_{mn}} \end{pmatrix},$$
 where $\delta_{ij}$ denotes the Kronecker delta symbol.

\newcommand\HC{{\mathbb H}\times{\mathbb C}}
\newcommand\ddx{{{\partial^2}\over{\partial x^2}}}
\newcommand\ddy{{{\partial^2}\over{\partial y^2}}}
\newcommand\ddu{{{\partial^2}\over{\partial u^2}}}
\newcommand\ddv{{{\partial^2}\over{\partial v^2}}}
\newcommand\px{{{\partial}\over{\partial x}}}
\newcommand\py{{{\partial}\over{\partial y}}}
\newcommand\pu{{{\partial}\over{\partial u}}}
\newcommand\pv{{{\partial}\over{\partial v}}}
\newcommand\pxu{{{\partial^2}\over{\partial x\partial u}}}
\newcommand\pyv{{{\partial^2}\over{\partial y\partial v}}}

\vskip 0.3cm C. L. Siegel [12] introduced the symplectic metric
$ds^2_n$ on ${\mathbb H}_n$ invariant under the action (1.1) of
$Sp(n,\BR)$ given by
\begin{eqnarray}
ds_n^2=\s \Big(Y^{-1}dZ\,Y^{-1}d{\overline Z}\Big)
\end{eqnarray}
and H. Maass [8] proved that the differential operator
\begin{eqnarray}
\Delta_n=4\,\s\left( Y \,\,^t\!\left( Y\PZB\right) \PZ\right)
\end{eqnarray}
is the Laplacian of ${\mathbb H}_n$ for the symplectic metric
$ds_n^2.$ Here $\sigma(A)$ denotes the trace of a square matrix
$A$.

\vskip 0.3cm In this paper, for arbitrary positive integers $n$
and $m$, we express the $G^J$-invariant metrics on $\SJ$ and their
Laplacians explicitly.

In fact, we prove the following theorems.

\vskip 0.5cm \noindent \begin{theorem}
For any two positive real numbers $A$ and $B$, the following
metric
\newcommand\zo{\overline{Z}}
\newcommand\zi{{Z^{-1}}}
\newcommand\yi{{Y^{-1}}}
\newcommand\zoi{{\,\overline{Z}^{-1}}}
\newcommand\tai{{(-\theta_2)^{-1}}}

\begin{eqnarray*}
ds_{n,m;A,B}^2&=& A\,\s\Big(\yi dZ\,\yi\bz\Big)\\
&&\ \ +\,B\,\bigg\{\s \Big(\yi\,{}^tV\,V\yi dZ\,\yi\bz\Big)+\,\s\Big(\yi\,{}^t(dW)\,\bw\Big)\\
&&\quad\quad-\s\Big(V\yi
dZ\,\yi\,{}^t(\bw)\Big)-\,\s\Big(V\yi\bz\,\yi\,{}^t(dW)\Big)\bigg\}
\end{eqnarray*}

\noindent is a Riemannian metric on $\Hnm$ which is invariant
under the action (1.2) of the Jacobi group $G^J$.
\end{theorem}

\vskip 0.3cm \noindent
\begin{theorem}
For any two positive real numbers $A$ and $B$, the Laplacian
$\Delta_{n,m;A,B}$ of $(\Hnm,ds^2_{n,m;A,B})$ is given by

\begin{eqnarray*}
\Delta_{n,m;A,B}&=&\frac4A\,\bigg\{ \s\left(\,Y\,\,
^t\!\left(Y\PZB\right)\PZ\,\right)\,+\,
\s\left(\,VY^{-1}\,^tV\ \,^t\!\left(Y\PWB\right)\,\PW\,\right)\\
& &\ \ +\,\s\left(V\
\,^t\!\left(Y\PZB\right)\PW\,\right)+\,\s\left(\,^tV \ \,
^t\!\left(Y\PWB\right)\PZ\,\right)
\bigg\}\\
&&\quad+\frac4B\,\, \s\left(Y\PW ^t\!\left(\PWB\right)\right).
\end{eqnarray*}

\newcommand\w{\wedge}

\noindent The following differential form
\begin{eqnarray*}
dv=\,\left(\,\text{det}\,Y\,\right)^{-(n+m+1)}[dX]\w [dY]\w [dU]\w
[dV]
\end{eqnarray*}
\noindent is a $G^J$-invariant volume element on $\Hnm$, where
$$[dX]=\w_{\mu\leq\nu}dx_{\mu\nu},\quad [dY]=\w_{\mu\leq\nu}
dy_{\mu\nu},\quad [dU]=\w_{k,l}du_{kl}\quad \text{and} \quad
[dV]=\w_{k,l}dv_{kl}.$$ The point is that the invariant metric
$ds_{n,m;A,B}^2$ and its Laplacian $\Delta_{n,m;A,B}$ are
expressed in terms of the {\it trace form}.
\end{theorem}

\vskip 0.2 cm For the case $n=m=1$ and $A=B=1$, Berndt proved in
[1]\,\,(cf.\,[19]) that the metric $ds^2_{1,1}$ on $\HC$ defined
by
\begin{eqnarray*} ds_{1,1}^2:=ds_{1,1;1,1}=\,&{{y\,+\,v^2}\over
{y^3}}\,(\,dx^2\,+\,dy^2\,)\,+\, {\frac 1y}\,(\,du^2\,+\,dv^2\,)\\
&\ \ -\,{{2v}\over {y^2}}\, (\,dx\,du\,+\,dy\,dv\,)
\end{eqnarray*}
is a Riemannian metric on $\HC$ invariant under the action (1.2)
of the Jacobi group and its Laplacian $\Delta_{1,1}$ is given by
\begin{eqnarray*} \Delta_{1,1}:=\Delta_{1,1;1,1}=\,& y^2\,\left(\,\ddx\,+\,\ddy\,\right)\,+\,
(\,y\,+\,v^2\,)\,\left(\,\ddu\,+\,\ddv\,\right)\\ &\ \
+\,2\,y\,v\,\left(\,\pxu\,+\,\pyv\,\right).
\end{eqnarray*}
\indent It is a pleasure to thank Eberhard Freitag for his helpful
advice and letting me know the paper [8] of H. Maass. \vskip 0.2cm
\noindent {\bf Notations:} \ \ We denote by $\BR$ and $\BC$ the
field of real numbers, and the field of complex numbers
respectively. The symbol ``:='' means that the expression on the
right is the definition of that on the left. For two positive
integers $k$ and $l$, $F^{(k,l)}$ denotes the set of all $k\times
l$ matrices with entries in a commutative ring $F$. For a square
matrix $A\in F^{(k,k)}$ of degree $k$, $\sigma(A)$ denotes the
trace of $A$. For any $M\in F^{(k,l)},\ ^t\!M$ denotes the
transpose matrix of $M$. $E_n$ denotes the identity matrix of
degree $n$. For $A\in F^{(k,l)}$ and $B\in F^{(k,k)}$, we set
$B[A]=\,^tABA.$ For a complex matrix $A$, ${\overline A}$ denotes
the complex {\it conjugate} of $A$. For $A\in \BC^{(k,l)}$ and
$B\in \BC^{(k,k)}$, we use the abbreviation $B\{
A\}=\,^t{\overline A}BA.$

\vskip 0.5cm
\begin{section}{{\bf \ Proof\ of\ Theorem\ 1.1}}
\setcounter{equation}{0}
\end{section}

Let $g=(M,(\l,\mu;\k))$ be an element of $G^J$ with
$M=\begin{pmatrix} A&B\\ C&D
\end{pmatrix}\in Sp(n,\BR)$ and $(Z,W)\in\Hnm$ with $Z\in\Hn$ and $W\in\Cmn.$
If we put $(Z_*,W_*):=g\cdot(Z,W),$ then we have
\begin{eqnarray*}
&Z_*=M\cdot Z=(AZ+B)(CZ+D)^{-1},\\
&W_*=(W+\l Z+\mu)(CZ+D)^{-1}.
\end{eqnarray*}
Thus we obtain
\begin{eqnarray}
dZ_*=dZ[(CZ+D)^{-1}]={}^t\!(CZ+D)^{-1}dZ(CZ+D)^{-1}
\end{eqnarray}
and
\begin{eqnarray}
\quad\quad\quad dW_*=dW(CZ+D)^{-1}+ \{ \l-(W+\l Z+\mu)(CZ+D)^{-1}C
\} dZ(CZ+D)^{-1}.
\end{eqnarray}
Here we used the following facts that
$$d(CZ+D)^{-1}=-(CZ+D)^{-1}C\,dZ(CZ+D)^{-1}$$
and that $(CZ+D)^{-1}C$ is symmetric.

\renewcommand\o{\overline}
We put
$$Z_*=X_*+iY_*,\quad W_*=U_*+iV_*,\quad X_*,Y_*,U_*,V_*\;\text{real}.$$
From [9], p.33 or [13], p.128, we know that
\begin{eqnarray}Y_*=Y\{(CZ+D)^{-1}\}={}^t\!(C\o{Z}+D)^{-1}\,Y(CZ+D)^{-1}.\end{eqnarray}
First of all, we recall that the following matrices
\begin{eqnarray*}
&&t(b)=\begin{pmatrix} E_n& b\\0&E_n\end{pmatrix},\quad b={}^tb\;\text{real},\\
&&g_0(h)=\begin{pmatrix} {}^th& 0\\0&h^{-1}\end{pmatrix},\quad h\in GL(n,\BR),\\
&&J_n=\begin{pmatrix} 0&-E_n\\E_n&0\end{pmatrix}
\end{eqnarray*}
generate the symplectic group $Sp(n,\mathbb R)$ (cf.\,[4],\,[5]).
Therefore the following elements\par\noindent
$t(b;\l,\mu,\k),\,g(h)$ and $\s_n$ of $G^J$ defined by
\begin{eqnarray*}
&& t(b;\l,\mu,\k)=\left(\begin{pmatrix} E_n&
b\\0&E_n\end{pmatrix},(\l,\mu;\k)
\right),\ \ b={}^tb\;\text{real},\;(\l,\mu;\k)\in H_{\BR}^{(n,m)},\\
&&g(h)=\left(\begin{pmatrix} {}^th&0\\0&h^{-1}\end{pmatrix},(0,0;0)\right),\ h\in GL(n,\BR),\\
&&\s_n=\left(\begin{pmatrix}
0&-E_n\\E_n&0\end{pmatrix},(0,0;0)\right)
\end{eqnarray*}
generate the Jacobi group $G^J.$ So it suffices to prove the
invariance of the metric $ds^2_{n,m;A,B}$ under the action of the
generators $t(b;\l,\mu,\k),\,g(h)$ and $\s_n.$ For brevity, we
write
\begin{eqnarray*}
&&(a)=\s\Big(Y^{-1}dZ\,Y^{-1}d\o{Z}\Big),\\
&&(b)=\s\Big(Y^{-1}{}^tVVY^{-1}dZ\,Y^{-1}d\o{Z}\Big),\\
&&(c)=\s\Big(Y^{-1}{}^t(dW)d\o{W}\Big),\\
&&(d)=-\s\Big(V\,Y^{-1}dZ\,Y^{-1}\,{}^t(d\o{W})\,+\,V\,Y^{-1}d\o{Z}\,Y^{-1}\,{}^t(dW)\Big)
\end{eqnarray*}
and
\begin{eqnarray*}
&&(a)_*=\s\Big(Y^{-1}_*dZ_*\,Y_*^{-1}d\o{Z}_*\Big),\\
&&(b)_*=\s\Big(Y_*^{-1}\,{}^tV_*V_*Y_*^{-1}\,dZ_*\,Y_*^{-1}d\o{Z}_*\Big),\\
&&(c)_*=\s\Big(Y_*^{-1}\,{}^t(dW_*)\,d\o{W}_*\Big),\\
&&(d)_*=-\s\Big(V_*\,Y_*^{-1}dZ_*\,Y_*^{-1}\,{}^t(d\o{W}_*)+\,V_*\,Y_*^{-1}d\o{Z}_*\,Y_*^{-1}\,{}^t(dW_*)\Big)
\end{eqnarray*}

\par
\vskip3mm {\bf Case I.} $g=t(b;\l,\mu,\k)$ with $b={}^tb$ real and
$(\l,\mu;\k)\in H_{\BR}^{(n,m)}.$ \vskip 0.2cm In this case, we
have
$$Z_*=Z+b,\quad Y_*=Y,\quad W_*=W+\l Z+\mu,\quad V_*=V+\l Y$$
and
$$dZ_* =dZ,\quad dW_*=dW+\l\, dZ.$$
Therefore
\begin{eqnarray*}
(a)_*&=&\s\Big(Y_*^{-1}dZ_*Y_*^{-1}d\o{Z_*}\Big)=\s\Big(Y^{-1}dZ\,Y^{-1}d\o{Z}\Big)=(a),\\
(b)_*&=&\s\Big(Y^{-1}\,{}^tV\,V\,Y^{-1}dZ\,Y^{-1}d\o{Z}\,\Big)+
\s\Big(Y^{-1}\,{}^tV\l\, dZ\,Y^{-1}d\o{Z}\,\Big)\\
&&\qquad+\s\Big(\,{}^t\l \,VY^{-1}dZ\,Y^{-1}d\o{Z}\,\Big)+\,\s\Big(\,{}^t\l\,\l\, dZ\,Y^{-1}d\o{Z}\,\Big),\\
(c)_*&=&\s \Big( Y^{-1}\,{}^t(dW)\,d\o{W}\, \Big)
+\s \Big( Y^{-1}\,{}^t(dW)\l\, d\o{Z} \,\Big)\\
&&\quad+\s \Big(Y^{-1}dZ\,\,{}^t\l\, d\o{W}\,\Big) +\s
\Big(Y^{-1}dZ\,\,{}^t\l\,\l\, d\o{Z}\,\Big)
\end{eqnarray*}
and
\begin{eqnarray*}
(d)_*&=&-\s\Big(V\,Y^{-1}dZ\,Y^{-1}\,{}^t(d\o{W}\,)\Big)-\s\Big(\l\,dZ \,Y^{-1}\,{}^t(d\o{W}\,)\Big)\\
&&\,-\s \Big( V\,Y^{-1}dZ\,Y^{-1}d\o{Z}\,\,{}^t\l \Big)
-\s \Big(\l\, dZ\,Y^{-1}d\o{Z}\,\,{}^t\l  \Big)\\
&&\,-\s\Big(V\, Y^{-1}d\o{Z}\,Y^{-1}\,{}^t(dW)\,\Big)-\s\Big(\l\,d\o{Z}\,Y^{-1}\,{}^t(dW)\Big)\\
&&\,-\s\Big(V\,Y^{-1}d\o{Z}\,Y^{-1}d{Z}\,\,{}^t\l\,\Big)
-\s\Big(\l\,d\o{Z}\,Y^{-1}d{Z}\,\,{}^t\l\Big).
\end{eqnarray*}
Thus we see that
\begin{eqnarray*}
(a)=(a)_*\quad \textrm{and}\quad (b)+(c)+(d)=(b)_*+(c)_*+(d)_*.
\end{eqnarray*}
Hence
\begin{eqnarray*}
ds_{n,m;A,B}^2=A\,(a)+B\Big\{ (b)+(c)+(d)\Big\}
\end{eqnarray*}

\noindent is invariant under the action of $t(B;\l,\mu,\k).$\par
\vskip 3mm {\bf Case II.} $g=g(h)$ with $h\in GL(n,\BR).$ \vskip
2mm In this case, we have
$$Z_*={}^th\,Z\,h,\quad Y_*={}^th\,Yh,\quad W_*=Wh,\quad V_*=Vh$$
and
$$dZ_*={}^th\,dZ\,h,\quad dW_*=dW\,h.$$
Therefore by an easy computation, we see that each of $(a),\ (b),\
(c)$ and $(d)$ is invariant under the action of all $g(h)$ with
$h\in GL(n,\BR).$ Hence the metric $ds_{n,m;A,B}^2$ is invariant
under the action of all $g(h)$ with $h\in GL(n,\BR).$
\par\vskip3mm
{\bf Case III.} $g=\s_n=\left(\begin{pmatrix}
0&-E_n\\E_n&0\end{pmatrix},(0,0;0)\right).$

\par
In this case, we have
\begin{equation}Z_*=-Z^{-1}\qquad\text{and}\qquad W_*=WZ^{-1}.\end{equation}
We set
$$\theta_1:=\text{Re}\,Z^{-1}\qquad\text{and}\qquad\theta_2:=\text{Im}\,Z^{-1}.$$
Then $\theta_1$ and $\theta_2$ are symmetric matrices and we have
\begin{equation}Y_*=-\theta_2\qquad\text{and}\qquad V_*:=\text{Im}\,W_*=V\theta_1+U\theta_2.\end{equation}
It is easy to see that
\begin{equation}Y=-Z\theta_2\o{Z}=-\o{Z}\theta_2Z,\end{equation}
\begin{equation}\theta_1Y+\theta_2X=0\end{equation}
and
\begin{equation}\theta_1X-\theta_2Y=E_n.\end{equation}
According to (2.6) and (2.7), we obtain
\begin{equation}X=(-\theta_2)^{-1}\theta_1Y\quad\quad \text{and}\quad\quad Y^{-1}=\theta_1(-\theta_2)^{-1}
\theta_1-\theta_2.\end{equation}
\newcommand\ta{\theta}
From (2.1) and (2.2), we have
\begin{equation}dZ_*=Z^{-1}dZZ^{-1}\end{equation}
and
\begin{equation}dW_*=dWZ^{-1}-WZ^{-1}dZZ^{-1}=\Big(dW-WZ^{-1}dZ\Big)Z^{-1}.\end{equation}
Therefore we have, according to (2.6) and (2.10),
\begin{eqnarray*}
(a)_*&=&\s \left( (-\th_2)^{-1}Z^{-1}dZ\,Z^{-1}(-\th_2)^{-1}\o{Z}^{-1}\bz\,\o{Z}^{-1} \right)\\
&=&\s\Big(Y^{-1}dZ\,Y^{-1}d\o{Z}\Big)=(a).
\end{eqnarray*}
According to (2.5)-(2.10), we have
\begin{eqnarray*}
(b)_* &=&\s \bigg(
(-\ta_2)^{-1}(\ta_1{}^tV+\ta_2{}^tU)(V\ta_1+U\ta_2)
(-\ta_2)^{-1}Z^{-1}dZ\,Y^{-1}d\o{Z}\,\o{Z}^{-1} \bigg)\\
&=&\s \bigg(
\{\,{}^tU-(-\ta_2)^{-1}\ta_1{}^tV\}\{U-V\ta_1(-\ta_2)^{-1}\}{Z}^{-1}
d{Z}\,Y^{-1}d\o{Z}\,\o{Z}^{-1} \bigg) \\
&=&\s \bigg( \{\,{}^t\o{W}+(iE_n-(-\ta_2)^{-1}\ta_1)\,{}^tV\}\{W-V(iE_n+\ta_1(-\ta_2)^{-1})\}\\
&&\qquad\quad Z^{-1}dZ\,Y^{-1}d\o{Z}\,\o{Z}^{-1} \bigg),
\end{eqnarray*}
\newcommand\zo{\overline{Z}}
\newcommand\zi{{Z^{-1}}}
\newcommand\yi{{Y^{-1}}}
\newcommand\zoi{{\,\overline{Z}^{-1}}}
\newcommand\tai{{(-\theta_2)^{-1}}}
\begin{eqnarray*}
(c)_* &=&\s \bigg(
(-\ta_2)^{-1}\Big(Z^{-1}\,{}^t(dW)-Z^{-1}dZZ^{-1}\,{}^tW\Big)\Big(d\o{W}\,
\o{Z}^{-1} -\o{W}\, \o{Z}^{-1}
d\o{Z}\,\o{Z}^{-1}\Big) \bigg)\\
&=&\s \bigg( (-\ta_2)^{-1}Z^{-1}\,{}^t(dW)\,\bw\zoi-\tai\zi\,{}^t(dW)\o{W}\zoi\bz\zoi\\
&&\qquad\quad-\tai Z^{-1}dZZ^{-1}\,{}^tW\bw\zoi\\
&&\qquad\quad+\tai Z^{-1}dZZ^{-1}\,{}^tW\o{W}\zoi\bz\zoi \bigg) \\
&=&\s \Big( Y^{-1}\,{}^t(dW)\,\bw \Big)-\s \Big( Y^{-1}\,{}^t(dW)\o{W}\zoi\bz \Big)\\
&&\qquad-\s \Big(Y^{-1}dZ\,Z^{-1}\,{}^tW\bw \,\Big) +\s
\Big(Y^{-1}dZ\,Z^{-1}\,{}^tW\o{W}\zoi\bz \Big)
\end{eqnarray*}
and
\begin{eqnarray*}
(d)_* &=&-\,\s\bigg( (V\ta_1+U\ta_2)\tai\zi dZ\zi\tai \left\{
\zoi\,\,{}^t(\bw\,)-\zoi\bz\zoi\,\,{}^t\o{W} \right\}
 \bigg)\\
&&\quad-\,\s \bigg((V\ta_1+U\ta_2) \tai\zoi\bz\zoi\tai \left\{
Z^{-1}\,{}^t(dW)-Z^{-1}dZZ^{-1}\,{}^tW
\right\} \bigg)\\
&=&-\,\s \bigg( (V\ta_1+U\ta_2)\tai\zi dZ\,\yi\,{}^t(\bw) \bigg)\\
&&\quad+\,\s \bigg((V\ta_1+U\ta_2) \tai\zi dZ\,\yi \bz\zoi\,{}^t\o{W} \bigg)\\
&&\quad-\,\s \bigg((V\ta_1+U\ta_2) \tai\zoi \bz\,\yi\,{}^t(dW)\bigg)\\
&&\quad+\,\s \bigg((V\ta_1+U\ta_2) \tai\zoi \bz\,\yi\, dZ
\zi\,{}^tW \bigg).
\end{eqnarray*}
\par
\newcommand\sq{\square}
\noindent Taking the $(dZ,\bw)$-part $\sq(Z,\o{W})$ in
$(b)_*+(c)_*+(d)_*,$ we have
\begin{eqnarray*}
\sq(Z,\o{W})
&=&-\s \Big(V\,\yi dZ\,\yi\,{}^t(\bw\,) \Big)+\s \Big(\yi dZ\,(\,{}^tW_*-\zi\, {}^tW)\,\bw \,\Big)\\
&=&-\s \Big(V\,\yi dZ\,\yi{}^t(\bw\,  ) \Big)\ \
\quad\text{because}\;W_*=W\zi\;(\text{cf.}\;(2.4)).
\end{eqnarray*}
Similiarly, if we take the $(\bz,dW)$-part $\sq(\o{Z},{W})$ in
$(b)_*+(c)_*+(d)_*,$ we have
\begin{eqnarray*}
\sq(\o{Z},{W})
&=&-\s \Big( V\,\yi \bz\,\yi\,{}^t(dW) \Big)+\s \Big( \bz\,\yi\,{}^t(dW)(\o{W_*}-\o{W}\zoi\,) \Big)\\
&=&-\s \Big( V\,\yi \bz\,\yi{}^t(dW) \Big)\ \
\quad\text{because}\;W_*=W\zi.
\end{eqnarray*}
If we take the  $(dW,\bw)$-part $\sq(W,\o{W})$ in
$(b)_*+(c)_*+(d)_*,$ we have
$$\sq (W,\o{W}\,)=\s \Big(\yi\,{}^t(dW)\,\bw \,\Big).$$
Finally, if we take the $(dZ,\bz)$-part $\sq(Z,\o{Z})$ in
$(b)_*+(c)_*+(d)_*,$ we have
\begin{eqnarray*}
\sq(Z,\o{Z})&=&\s \bigg(
\left\{\,{}^t\o{W}+(iE_n-\tai\ta_1)\,{}^tV \right\}
\left\{ W-V(iE_n+\ta_1\tai\,) \right\}\\
&&\quad\quad\zi dZ\,\yi\bz\zoi \bigg)\\
&&\quad+\,\s \left(\zi\,{}^tW\o{W}\zoi\bz\,\yi dZ \right)\\
&&\quad+\,\s \bigg(  {}^t\o{W}(V\ta_1+U\ta_2)\tai\zi dZ\,\yi\bz\zoi \bigg)\\
&&\quad+\,\s \bigg( {}^t{W}(V\ta_1+U\ta_2)\tai\zoi \bz\,\yi dZ\zi
\bigg).
\end{eqnarray*}
Since
\begin{eqnarray*}
(V\ta_1+U\ta_2)\tai&=&-U+V\ta_1\tai\\
&=&-W+V\{iE_n+\ta_1\tai\}\\
&=&-\o{W}-V\{iE_n-\ta_1\tai\},
\end{eqnarray*}
we have
\begin{eqnarray*}
\sq(Z,\o{Z})
&=&\s \bigg( \zoi \left\{ iE_n-\tai\ta_1 \right\}{}^tV \left\{W-V(iE_n+\ta_1\tai) \right\}\\
&&\quad\quad\zi dZ\,\yi\bz \bigg)\\
&&\quad-\s \bigg( \zoi \left\{ iE_n-\tai\ta_1 \right\}{}^tVW\zi dZ\,\yi\bz \bigg) \\
&=&-\s \bigg(  \zoi \left\{iE_n-\tai\ta_1 \right\}{}^tVV \left\{iE_n+\ta_1\tai \right\}\\
&&\quad\quad\zi dZ\,\yi\bz \bigg).
\end{eqnarray*}
By the way, according to (2.9), we obtain
\begin{eqnarray*}
\zoi\Big\{iE_n-\tai\ta_1\Big\}&=&(\ta_1-i\ta_2)\Big\{iE_n-\tai\ta_1\Big\}\\
&=&\ta_2-\ta_1\tai\ta_1=-\yi
\end{eqnarray*}
and
\begin{eqnarray*}
\Big\{iE_n+\ta_1\tai\Big\}\zi =\ta_1\tai\ta_1-\ta_2=\yi.
\end{eqnarray*}
Therefore
$$\sq(Z,\o{Z})=\s\Big(\yi\,{}^tV\,V\,\yi dZ\,\yi\,\bz\Big).$$
Hence $(a)=(a)_*$ and
\begin{eqnarray*}
(b)_*+(c)_*+(d)_*&=&\,\sq(Z,\o{W})+\sq(\o{Z},W)+\sq(W,\o{W})
+\sq(Z,\o{Z})\\
&=&\,(b)+(c)+(d).
\end{eqnarray*}
This implies that the metric
\begin{eqnarray*}
ds_{n,m;A,B}^2\,=A\, (a)\,+\,B\,\Big\{ (b)+(c)+(d)\Big\}
\end{eqnarray*}

\noindent is invariant under the action (1.2) of $\s_n.$\par
Consequently $ds_{n,m;A,B}^2$ is invariant under the action (1.2)
of the Jacobi group $G^J.$ In particular, for $(Z,W)=(iE_n,0),$ we
have
\begin{eqnarray*}
ds_{n,m;A,B}^2&=&\,A\cdot \s\Big(dZ\,\bz\Big)\,+\,B\cdot \s \Big( {}^t(dW)\bw \Big)\\
&=&\,A\,\Bigg\{\sum_{\mu=1}^n(dx_{\mu\mu}^2+dy_{\mu\mu}^2)+2\sum_{1\leq\mu<\nu\leq
n}
(dx_{\mu\nu}^2+dy_{\mu\nu}^2)\,\Bigg\}\\
&&\quad+\,B\,\Bigg\{\sum_{ 1\leq k\leq m,\, 1\leq l\leq n}
(du_{kl}^2+dv_{kl}^2)\Bigg\},
\end{eqnarray*}
which is clearly positive definite. Since $G^J$ acts on $\mathbb
H_{n,m}$ transitively, $ds_{n,m;A,B}^2$ is positive definite
everywhere in $\mathbb H_{n,m}.$ This completes the proof of
Theorem 1.1.\hfill $\square$

\vskip 0.3cm
\par\noindent
{\bf Remark\ 2.1.} The scalar curvature of the Siegel-Jacobi space
$(\BH_{n,m},\,ds_{n,m;A,B}^2)$ is constant because of the
transitive group action of $G^J$ on $\BH_{n,m}.$ In the special
case $n=m=1$ and $A=B=1$, by a direct computation, we see that the
scalar curvature of
$(\BH_{1,1}, ds_{1,1;1,1}^2)$ is $-3.$\\

\vskip 0.53cm
\begin{section}{{\bf \ Proof\ of\ Theorem\ 1.2}}
\setcounter{equation}{0}
\end{section}

If $(Z_*,W_*)= g \cdot (Z,W)$ with $g=\left( \left(\begin{matrix}
A&B\\ C&D\end{matrix}\right),(\l,\mu;\kappa)\right)\in G^J,$ we
can see easily that
\begin{eqnarray}
\frac{\partial}{\partial Z_*}&=&(CZ+D)\,\, {}^t\!\left\{
(CZ+D)\,\PZ\right\} \\
&&\quad +\,(CZ+D)\,\,{}^t\!\left\{
\Big(C\,{}^tW+C\,\,{}^t\mu-D\,\,{}^t\l\Big)\,\,{}^t\!\left(\PW\right)\right\}\nonumber
\end{eqnarray}
and
\begin{equation}\frac{\partial}{\partial W_*}=(CZ+D)\PW.\end{equation}
For brevity, we put
\begin{eqnarray*}
&&(\alpha):=4\,\s\left(\,Y\,\,
^t\!\left(Y\PZB\right)\PZ\,\right),\\
&&(\beta):=\;4\,\s\left(\, Y\,\PW\,\,{}^t\!\left( \PWB\right)\,\right),\\
&&(\g):=\;4\,\s\left(\,VY^{-1}\,^tV\,\,^t\!\left(Y\PWB\right)\,\PW\,\right),\\
&&(\delta):=\;4\,\s\left(V\,\,^t\!\left(Y\PZB\right)\PW\,\right)
\end{eqnarray*}
and
\begin{eqnarray*}
(\epsilon):=\;4\,\s\left(\,^tV\,\,^t\!\left(Y\PWB\right)\PZ\,\right).
\end{eqnarray*}
We also set
\begin{eqnarray*}
&&(\alpha)_*:=4\,\s\left(\,Y_*\,\,
^t\!\left(Y_*\PZB_*\right)\PZ_*\,\right),\\
&&(\beta)_*:=\;4\,\s\left(\, Y_*\,\PW_*\,\,{}^t\!\left( \PWB_*\right)\,\right),\\
&&(\g)_*:=\;4\,\s\left(\,V_*\,Y_*^{-1}\,^tV_*\,\,^t\!\left(Y_*\PWB_*\right)\,\PW_*\,\right),\\
&&(\delta)_*:=\;4\,\s\left(V_*\,\,^t\!\left(Y_*\PZB_*\right)\PW_*\,\right)
\end{eqnarray*}
and
\begin{eqnarray*}
(\epsilon)_*:=\;4\,\s\left(\,^tV_*\,\,^t\!\left(Y_*\PWB_*\right)\PZ_*\,\right).
\end{eqnarray*}

We need the following lemma for the proof of Theorem 1.2. H.
Maass\,[8] observed the following useful fact. \vskip 0.3cm
\noindent
\begin{lemma}
(a) Let $A$ be an $n\times k$ matrix and let $B$ be a $k\times n$
matrix. Assume that the entries of $A$ commute with the entries of
$B$. Then $\s (AB)=\s (BA).$ \vskip 0.1cm\noindent (b) Let $A$ be
an $m\times n$ matrix and $B$ an $n\times l$ matrix. Assume that
the entries of $A$ commute with the entries of $B$. Then
${}^t(AB)=\,{ }^tB\,\,{ }^tA.$ \vskip 0.1cm \noindent (c) Let
$A,\,B$ and $C$ be a $k\times l$, an $n\times m$ and an $m\times
l$ matrix respectively. Assume that the entries of $A$ commute
with the entries of $B$. Then
\begin{eqnarray*}
{ }^t(A\,\,{ }^t(BC))=\,B\,\,{ }^t(A\,\,^tC).
\end{eqnarray*}
\end{lemma}

\noindent {\it Proof.} The proof follows immediately from the
direct computation. \hfill$\square$ \vskip 0.3cm \indent Now we
are ready to prove Theorem 1.2. First of all, we shall prove that
$\Delta_{n,m;A,B}$ is invariant under the action of the generators
$t(b;\l,\mu,\k), \;g(h)$ and $\s_n.$ \vskip 0.2cm \noindent
{\bf Case\ I.} $g=t(b;\l,\mu,\k) =\left( \left(\begin{matrix} E_n&b\\
0&E_n\end{matrix}\right),(\l,\mu;\k)\right)$ with $b=\,^tb$
real.\par In this case, we have
$$Y_*=Y,\qquad V_*=V+\l\,Y$$
and
$$ \frac{\partial}{\partial Z_*}=\PZ-\,\, ^t\!\left(\,{ }^t\l\,\, ^t\!\left(\PW\right)\right) \quad\text{and}\quad \;\;
\frac{\partial}{\partial W_*}=\PW\,.$$ Using Lemma 3.1, we obtain
\begin{eqnarray*}
(\alpha)_*&=&(\alpha)-\s\left( \l\, Y\,\,^t\!\left(Y\PZB\right)\PW       \right)\\
& & \ -\, \s\left(Y\,\,^t\!\l\,\,^t\!\!\left( Y\PWB\right)\PZ  \right)\,+\,\s\left(\l\,Y\,\,^t\!\l\,\,^t\!\!\left( Y\PWB\right)\PW  \right)\\
(\beta)_*&=& \, \s\left(Y_*\frac{\partial}{\partial W_*}
^t\left(\frac{\partial}{\partial \overline {W_*}}\right)\right)=(\beta),\\
(\g)_*&=&(\g)\,+\, \s \left( \l\,\,^tV\,\,^t\!\!\left( Y\PWB\right)\PW \right)\\
& & \ +\,\s\left( V\,\,^t\!\l\,\,^t\!\!\left(
Y\PWB\right)\PW\right)\,+\,\s \left( \l\,
Y\,\,^t\!\l\,\,^t\!\!\left( Y\PWB\right)\PW\right),\\
(\delta)_*&=&(\delta)+\,\s\left( \l\, Y\,\,^t\!\left(Y\PZB\right)\PW       \right)\\
& &\ -\,\s\left( V\,\,^t\!\l\,\,^t\!\!\left(
Y\PWB\right)\PW\right)\,+\,\s \left( \l\,
Y\,\,^t\!\l\,\,^t\!\!\left( Y\PWB\right)\PW\right)
\end{eqnarray*}
and
\begin{eqnarray*}
(\epsilon)_*&=&(\epsilon)+\,\s\left(Y\,\,^t\!\l\,\,^t\!\!\left(
Y\PWB\right)\PZ  \right)\\
& &\ -\,\s \left( \l\,\,^tV\,\,^t\!\!\left( Y\PWB\right)\PW
\right)\,-\,\s \left( \l\, Y\,\,^t\!\l\,\,^t\!\!\left(
Y\PWB\right)\PW\right).
\end{eqnarray*}
Thus $(\beta)=(\beta)_*$ and
\begin{eqnarray*}
(\alpha)+(\g)+(\delta)+(\epsilon)=(\alpha)_*
+(\g)_*+(\delta)_*+(\epsilon)_*.
\end{eqnarray*}

\noindent Hence
\begin{eqnarray*}
\Delta_{n,m;A,B}= {\frac 4B}\,(\beta)\,+\,{\frac
4A}\,\Big\{(\alpha)+(\g)+(\delta)+(\epsilon)\Big\}
\end{eqnarray*}

\noindent is invariant under the action of all $\,t(b;\l,\mu,\k).$

\vskip 0.3cm \noindent
{\bf Case\ II.} $g=g(h)=\left( \left(\begin{matrix} ^th&0\\
0&h^{-1}\end{matrix}\right),(0,0;0)\right)$ with $h\in
GL(n,\BR)$.\par In this case, we have
$$Y_*=\,^th\,Y\,h\quad V_*=V\,h$$
and
$$ \frac{\partial}{\partial Z_*}=h^{-1}\,\,^t\!\left(h^{-1}\PZ\right),\quad\quad \;\;
\frac{\partial}{\partial W_*}=h^{-1}\PW.$$ According to Lemma 3.1,
we see that each of $(\alpha),\,(\beta),\,(\g),\,(\delta)$ and
$(\epsilon)$ is invariant under the action of all $g(h)$ with
$h\in GL(n,\BR).$ Therefore $\Delta_{n,m;A,B}$ is invariant under
the action of all $g(h)$ with $h\in GL(n,\BR).$

\vskip 0.3cm \noindent
{\bf Case\ III.} $g=\s_n=\left( \left(\begin{matrix} 0& -E_n\\
E_n&0\end{matrix}\right),(0,0;0)\right).$\par In this case, we
have
$$Z_*=-Z^{-1} \;\;\quad \text{and} \quad\;\;W_*=WZ^{-1}.$$
We set
$$ \th_1:= \text{Re} \,Z^{-1}\;\;\quad\text{and}\quad\;\; \th_2:= \text{Im} \, Z^{-1}.$$
Then we obtain the relations (2.5)-(2.9). From (2.6), we have the
relation
\begin{equation} \theta_2{\overline Z}=-Z^{-1}Y.\end{equation}
It follows from the relation (2.3) that
\begin{equation} Y_*={\overline Z}^{-1}YZ^{-1}=Z^{-1}Y{\overline
Z}^{-1}=-\theta_2.\end{equation} From (2.9), we obtain
\begin{equation} \theta_1\theta_2^{-1}\theta_1=-Y^{-1}-\theta_2.\end{equation}
According to (3.1) and (3.2), we have
\begin{equation} \frac{\partial}{\partial Z_*}=Z\,\,^t\!\!\left( Z\PZ\right)\,
+\,Z\,\,^t\!\!\left(\,^tW\,\,^t\!\!\left(\PW\right)\right)
\end{equation} and
\begin{equation} \frac{\partial}{\partial W_*}=Z\PW.\end{equation}
\par
From (2.6), (3.3) and Lemma 3.1, we obtain
\begin{eqnarray*}
(\alpha)_*&=& (\alpha)\,-\,\s\bigg( U\theta_2{\overline
Z}\,\,^t\!\!\left(Y\PZB\right)\PW \bigg) \,-\,i\,\s\bigg(
V\theta_2{\overline
Z}\,\,^t\!\!\left(Y\PZB\right)\PW\bigg)\\
& &\ -\s\bigg(
Z\theta_2\,\,^tU\,\,^t\!\!\left(Y\PWB\right)\PZ\bigg)\,+\,i\,\s\bigg(
Z\theta_2\,\,^tV\,\,^t\!\!\left(Y\PWB\right)\PZ\bigg)\\
& &\ -\,\s\bigg(
W\theta_2\,\,^tU\,\,^t\!\!\left(Y\PWB\right)\PW\bigg)\,+\,\s\bigg(
W\theta_2\,\,^tV\,\,^t\!\!\left(Y\PWB\right)\PW\bigg).
\end{eqnarray*}
From the relation (3.4), we see $(\beta)_*=(\beta)$. According to
(3.3),\,(3.5) and Lemma 3.1, we otain
\begin{eqnarray*}
(\g)_*=(\g)\,+\,\s\bigg(
(V\theta_2\,^tV-V\theta_1\,^tU-U\theta_1\,^tV-U\theta_2\,^tU)\,\,^t\!\!\left(Y\PWB\right)\PW\bigg).
\end{eqnarray*}
Using the relation (3.3) and Lemma 3.1, we finally obtain
\begin{eqnarray*}
(\delta)_*&=&\,\s\bigg( V\theta_1{\overline
Z}\,\,^t\!\!\left(Y\PZB\right)\PW \bigg) \,+\,\s\bigg(
U\theta_2{\overline
Z}\,\,^t\!\!\left(Y\PZB\right)\PW\bigg)\\
& &\ +\,\s\bigg( V\theta_1\,\,^t{\overline
W}\,\,^t\!\!\left(Y\PWB\right)\PW\bigg)\,+\,\s\bigg(
U\theta_2\,^t{\overline W}\,\,^t\!\!\left(Y\PWB\right)\PW\bigg)
\end{eqnarray*}
and
\begin{eqnarray*}
(\epsilon)_*&=&\,\s\bigg(
Z\theta_1\,^tV\,\,^t\!\!\left(Y\PWB\right)\PZ \bigg) \,+\,\s\bigg(
Z\theta_2\,^tU\,\,^t\!\!\left(Y\PWB\right)\PZ\bigg)\\
& &\ +\,\s\bigg(
W\theta_1\,\,^tV\,\,^t\!\!\left(Y\PWB\right)\PW\bigg)\,+\,\s\bigg(
W\theta_2\,^tU\,\,^t\!\!\left(Y\PWB\right)\PW\bigg).
\end{eqnarray*}

Using the fact $Z^{-1}=\theta_1+i\theta_2$, we can show that
$$ (\alpha)+(\g)+(\delta)+(\epsilon)\,=\,(\alpha)_*+(\g)_*+(\delta)_*+(\epsilon)_*.$$

\noindent Hence
\begin{eqnarray*}
\Delta_{n,m;A,B}= {\frac 4B}\,(\beta)\,+\,{\frac
4A}\,\Big\{(\alpha)+(\g)+(\delta)+(\epsilon)\Big\}
\end{eqnarray*}

\noindent is invariant under the action of $\s_n$.\par
Consequently $\Delta_{n,m;A,B}$ is invariant under the action
(1.2) of $G^J$. In particular, for $(Z,W)=(iE_n,0),$ the
differential operator $\Delta_{n,m;A,B}$ coincides with the
Laplacian for the metric $ds^2_{n,m;A,B}$. It follows from the
invariance of $\Delta_{n,m;A,B}$ under the action (1.2) and the
transitivity of the action of $G^J$ on $\Hnm$ that
$\Delta_{n,m;A,B}$ is the Laplacian of $(\Hnm,\,ds_{n,m;A,B}^2).$
The invariance of the differential form $dv$ follows from the fact
that the following differential form
$$(\,\det Y\,)^{-(n+1)}[dX]\wedge [dY]$$
is invariant under the action (1.1) of $Sp(n,\BR)$\,(cf.\,[13],
p.\,130). \hfill$\square$

\vskip 0.73cm
\begin{section}{{\bf  Remark on Spectral Theory of $\Delta_{n,m;A,B}$ on Siegel-Jacobi Space }}
\setcounter{equation}{0}
\end{section}
\newcommand\CCF{{\mathcal F}_n}
\newcommand\CM{{\mathcal M}}
\newcommand\CHX{{\mathcal H}_{\xi}}
\newcommand\Ggh{\Gamma_{g,h}}
\newcommand\CP{{\mathcal P}_n}
\newcommand\Dnm{{\mathbb D}_n \times {\mathbb C}^{(m,n)}}
\newcommand\BHn{\BH_n}
\newcommand\BDn{\BD_n}
\newcommand\Cal{\mathcal}
\vskip 0.3cm Before we describe a fundamental domain for the
Siegel-Jacobi space, we review the Siegel's fundamental domain for
the Siegel upper half plane. \vskip 0.2cm We let
\begin{equation*}
{\mathcal P}_n=\left\{ Y\in\BR^{(n,n)}\,|\ Y=\,^tY>0\,\right\}
\end{equation*}

\noindent be an open cone in $\BR^{n(n+1)/2}$. The general linear
group $GL(n,\BR)$ acts on ${\mathcal P}_n$ transitively by
\begin{equation*}
h\circ Y:=h\,Y\,{}^th,\quad h\in GL(n,\BR),\ Y\in {\mathcal P}_n.
\end{equation*}

\noindent Thus ${\mathcal P}_n$ is a symmetric space diffeomorphic
to $GL(n,\BR)/O(n).$ We let
\begin{eqnarray*}
GL(n,{\mathbb Z})=\Big\{ h\in GL(n,\BR)\,\Big| \ h\ \textrm{is
integral}\ \Big\}
\end{eqnarray*}

\noindent be the discrete subgroup of $GL(n,\BR)$.

\newcommand\Mg{{\mathcal M}_n}
\newcommand\Rg{{\mathcal R}_n}
\newcommand\ba{\backslash}
\newcommand\BZ{\mathbb Z}
\newcommand\Om{\Omega}
\vskip 0.2cm The fundamental domain $\Rg$ for $GL(n,\BZ)\ba
{\mathcal P}_n$ which was found by H. Minkowski\,[10] is defined
as a subset of ${\mathcal P}_n$ consisting of $Y=(y_{ij})\in
{\mathcal P}_n$ satisfying the following conditions (M.1)-(M.2)\
(cf.\,[9] p.\,123): \vskip 0.1cm (M.1)\ \ \ $aY\,^ta\geq y_{kk}$\
\ for every $a=(a_i)\in\BZ^n$ in which $a_k,\cdots,a_n$ are
relatively prime for\par \ \ \ \ \ \ \ \ \ \ $k=1,2,\cdots,n$.
\vskip 0.1cm (M.2)\ \ \ \ $y_{k,k+1}\geq 0$ \ for
$k=1,\cdots,n-1.$ \vskip 0.1cm  \noindent We say that a point of
$\Rg$ is {\it Minkowski reduced} or simply {\it M}-{\it reduced}.

\vskip 0.1cm Siegel\,[12] determined a fundamental domain
${\mathcal F}_n$ for $\G_n\ba \BH_n,$ where $\G_n=Sp(n,\BZ)$ is
the Siegel modular group of degree $n$. We say that $\Om=X+iY\in
\BH_n$ with $X,\,Y$ real is {\it Siegel reduced} or {\it S}-{\it
reduced} if it has the following three properties: \vskip 0.1cm
(S.1)\ \ \ $\det (\text{Im}\,(\g\cdot\Om))\leq \det
(\text{Im}\,(\Om))\qquad\text{for\ all}\ \g\in\G_n$; \vskip 0.1cm
(S.2)\ \ $Y=\text{Im}\,\Om$ is M-reduced, that is, $Y\in \Rg\,;$
\vskip 0.1cm (S.3) \ \ $|x_{ij}|\leq {\frac 12}\quad \text{for}\
1\leq i,j\leq n,\ \text{where}\ X=(x_{ij}).$ \vskip 0.21cm
${\mathcal F}_n$ is defined as the set of all Siegel reduced
points in $\BH_n.$ Using the highest point method, Siegel [12]
proved the following (F1)-(F3)\,(cf. [9], p.\,169): \vskip 0.1cm
(F1)\ \ \ $\G_n\cdot {\mathcal F}_n=\BH_n,$ i.e.,
$\BH_n=\cup_{\g\in\G_n}\g\cdot {\mathcal F}_n.$ \vskip 0.1cm (F2)\
\ ${\mathcal F}_n$ is closed in $\BH_n.$ \vskip 0.1cm (F3)\ \
${\mathcal F}_n$ is connected and the boundary of ${\mathcal F}_n$
consists of a finite number of hyperplanes.

\vskip 0.2cm The metric $ds_n^2$ given by (1.3) induces a metric
$ds_{{\mathcal F}_n}^2$ on ${\mathcal F}_n$. Siegel\,[12] computed
the volume of ${\mathcal F}_n$
\begin{equation*}
\text{vol}\,(\CCF)=2\prod_{k=1}^n\pi^{-k}\G
(k)\zeta(2k),\end{equation*} where $\G (s)$ denotes the Gamma
function and $\zeta (s)$ denotes the Riemann zeta function. For
instance,
$$\text{vol}\,({\mathcal F}_1)={{\pi}\over 3},\quad \text{vol}\,({\mathcal F}_2)={{\pi^3}\over {270}},
\quad \text{vol}\,({\mathcal F}_3)={{\pi^6}\over {127575}},\quad
\text{vol}\,({\mathcal F}_4)={{\pi^{10}}\over {200930625}}.$$

\vskip 0.2cm Let $f_{kl}\,(1\leq k\leq m,\ 1\leq l\leq n)$ be the
$m\times n$ matrix with entry $1$ where the $k$-th row and the
$l$-th column meet, and all other entries $0$. For an element
$\Om\in \BH_n$, we set for brevity
\begin{equation*}
h_{kl}(\Om):=f_{kl}\Om,\qquad 1\leq k\leq m,\ 1\leq l\leq
n.\end{equation*}

For each $\Om\in {\mathcal F}_n,$ we define a subset $P_{\Om}$ of
$\BC^{(m,n)}$ by
\begin{equation*}
P_{\Om}=\left\{ \,\sum_{k=1}^m\sum_{j=1}^n \l_{kl}f_{kl}+
\sum_{k=1}^m\sum_{j=1}^n \mu_{kl}h_{kl}(\Om)\,\Big|\ 0\leq
\l_{kl},\mu_{kl}\leq 1\,\right\}. \end{equation*}

\noindent For each $\Om\in {\mathcal F}_n,$ we define the subset
$D_{\Om}$ of $\BH_n\times \BC^{(m,n)}$ by
\begin{equation*} D_{\Om}:=\left\{\,(\Om,Z)\in\BH_n\times \BC^{(m,n)}\,\vert\ Z\in
P_{\Om}\,\right\}.\end{equation*}
\newcommand\Fgh{{\mathcal F}_{n,m}}
\noindent We define
\begin{equation*} \Fgh:=\cup_{\Om\in {\mathcal F}_n}D_{\Omega}.\end{equation*}

\vskip 0.3cm
\begin{theorem}
Let
$$\Gamma_{n,m}:=Sp(n,{\mathbb Z})\ltimes H_{\mathbb Z}^{(n,m)}$$
be the discrete subgroup of $G^J$, where
$$H_{\BZ}^{(n,m)}=\left\{ (\lambda,\mu;\kappa)\in
H_{\BR}^{(n,m)}\,|\ \lambda,\mu,\kappa \ \textrm{are integral}\
\right\}.$$ Then $\Fgh$ is a fundamental domain for $\G_{n,m}\ba
\BH_{n,m}.$
\end{theorem}
\vskip 0.2cm \noindent $\textit{Proof.}$ The proof can be found in
[20]. \hfill $\square$

\vskip 0.3cm In the case $n=m=1$, R. Berndt [2] introduced the
notion of Maass-Jacobi forms. Now we generalize this notion to the
general case.

\vskip 0.2cm
\begin{definition}
For brevity, we set $\Delta_{n,m}:=\Delta_{n,m;1,1}$ (cf. Theorem
1.2). Let
$$\Gamma_{n,m}:=Sp(n,{\mathbb Z})\ltimes H_{\mathbb Z}^{(n,m)}$$
be the discrete subgroup of $G^J$, where
$$H_{\BZ}^{(n,m)}=\left\{ (\lambda,\mu;\kappa)\in
H_{\BR}^{(n,m)}\,|\ \lambda,\mu,\kappa \ \textrm{are integral}\
\right\}.$$ A smooth function $f:\Hnm\lrt \BC$ is called a
$\textsf{Maass}$-$\textsf{Jacobi form}$ on $\Hnm$ if $f$ satisfies
the following conditions (MJ1)-(MJ3)\,:\vskip 0.1cm (MJ1)\ \ \ $f$
is invariant under $\G_{n,m}.$\par (MJ2)\ \ \ $f$ is an
eigenfunction of the Laplacian
 $\Delta_{n,m}$.\par (MJ3)\ \ \ $f$
has a polynomial growth, that is, there exist a constant $C>0$ and
a positive \par \ \ \ \ \ \ \ \ \ \ \ integer $N$ such that
\begin{equation*}
|f(X+iY,Z)|\leq C\,|p(Y)|^N\quad \textrm{as}\ \det
Y\longrightarrow \infty,
\end{equation*}

\ \ \ \ \ \ \ \ \ \ \ where $p(Y)$ is a polynomial in
$Y=(y_{ij}).$
\end{definition}

\vskip 0.3cm It is natural to propose the following problems.

\vskip 0.3cm\noindent {\bf {Problem\ A}\,:} Construct Maass-Jacobi
forms.

\vskip 0.3cm\noindent {\bf {Problem\ B}\,:} Find all the
eigenfunctions of $\Delta_{n,m}.$

\vskip 0.3cm  We consider the simple case $n=m=1.$ A metric
$ds_{1,1}^2$ on ${\mathbf H}_1\times \BC$ given by
\begin{align*} ds^2_{1,1}\,=\,&{{y\,+\,v^2}\over
{y^3}}\,(\,dx^2\,+\,dy^2\,)\,+\, {\frac 1y}\,(\,du^2\,+\,dv^2\,)\\
&\ \ -\,{{2v}\over {y^2}}\, (\,dx\,du\,+\,dy\,dv\,)\end{align*} is
a $G^J$-invariant K{\"a}hler metric on ${\mathbf H}_1\times \BC$.
Its Laplacian $\Delta_{1,1}$ is given by
\begin{align*} \Delta_{1,1}\,=\,& \,y^2\,\left(\,\ddx\,+\,\ddy\,\right)\,\\
&+\, (\,y\,+\,v^2\,)\,\left(\,\ddu\,+\,\ddv\,\right)\\ &\ \
+\,2\,y\,v\,\left(\,\pxu\,+\,\pyv\,\right).\end{align*}

\vskip 0.2cm We provide some examples of eigenfunctions of
$\Delta_{1,1}$. \vskip 0.2cm (1) $h(x,y)=y^{1\over
2}K_{s-{\frac12}}(2\pi |a|y)\,e^{2\pi iax} \ (s\in \BC,$
$a\not=0\,)$ with eigenvalue $s(s-1).$ Here
$$K_s(z):={\frac12}\int^{\infty}_0 \exp\left\{-{z\over
2}(t+t^{-1})\right\}\,t^{s-1}\,dt,$$ \indent \ \ \ where
$\mathrm{Re}\,z
> 0.$ \par (2) $y^s,\ y^s x,\ y^s u\ (s\in\BC)$ with eigenvalue
$s(s-1).$
\par
 (3) $y^s v,\ y^s uv,\ y^s xv$ with eigenvalue $s(s+1).$
\par
(4) $x,\,y,\,u,\,v,\,xv,\,uv$ with eigenvalue $0$.
\par
(5) All Maass wave forms.

\vskip 0.3cm We fix two positive integers $m$ and $n$ throughout
this section. \vskip 0.1cm For an element $\Om\in \BH_n,$ we set
\begin{equation*}
L_{\Om}:=\BZ^{(m,n)}+\BZ^{(m,n)}\Om\end{equation*} \noindent It
follows from the positivity of $\text{Im}\,\Om$ that the elements
$f_{kl},\,h_{kl}(\Om)\,(1\leq k\leq m,\ 1\leq l\leq n)$ of
$L_{\Om}$ are linearly independent over $\BR$. Therefore $L_{\Om}$
is a lattice in $\BC^{(m,n)}$ and the set
$\left\{\,f_{kl},\,h_{kl}(\Om)\,|\ 1\leq k\leq m,\ 1\leq l\leq n\,
\right\}$ forms an integral basis of $L_{\Om}$. We see easily that
if $\Om$ is an element of $\BH_n$, the period matrix
$\Om_*:=(I_n,\Om)$ satisfies the Riemann conditions (RC.1) and
(RC.2)\,: \vskip 0.1cm (RC.1) \ \ \ $\Om_*J_n\,^t\Om_*=0\,$;
\vskip 0.1cm (RC.2) \ \ \ $-{1 \over
{i}}\Om_*J_n\,^t{\overline{\Om}}_*
>0$.

\vskip 0.2cm \noindent Thus the complex torus
$A_{\Om}:=\BC^{(m,n)}/L_{\Omega}$ is an abelian variety.  \vskip
0.2cm It might be interesting to investigate the spectral theory
of the Laplacian $\Delta_{n,m}$ on a fundamental domain ${\mathcal
F}_{n,m}$. But this work is very complicated and difficult at this
moment. It may be that the first step is to develop the spectral
theory of the Laplacian $\Delta_{\Omega}$ on the abelian variety
$A_{\Omega}.$ The second step will be to study the spectral theory
of the Laplacian $\Delta_n$\,(see (1.4)) on the moduli space
$\Gamma_n\backslash \BH_n$ of principally polarized abelian
varieties of dimension $n$. The final step would be to combine the
above steps and more works to develop the spectral theory of the
Lapalcian $\Delta_{n,m}$ on ${\mathcal F}_{n,m}.$ $
\textit{Maass-Jacobi forms}$ play an important role in the
spectral theory of $\Delta_{n,m}$ on ${\mathcal F}_{n,m}.$ Here we
deal only with the spectral theory $\Delta_{\Omega}$ on
$L^2(A_{\Omega}).$

 \vskip 0.1cm We fix
an element $\Om=X+iY$ of $\BH_n$ with $X=\text{Re}\,\Om$ and
$Y=\text{Im}\, \Om.$ For a pair $(A,B)$ with $A,B\in\BZ^{(m,n)},$
we define the function $E_{\Om;A,B}:\BC^{(m,n)}\lrt \BC$ by
\begin{equation*}
E_{\Om;A,B}(Z)=e^{2\pi i\left( \s(\,^tAU\,)+\, \s
((B-AX)Y^{-1}\,^tV)\right)},\end{equation*} where $Z=U+iV$ is a
variable in $\BC^{(m,n)}$ with real $U,V$. \vskip 0.1cm\noindent
\begin{lemma} For any $A,B\in \BZ^{(m,n)},$ the function
$E_{\Om;A,B}$ satisfies the following functional equation
\begin{equation*}
E_{\Om;A,B}(Z+\l \Om+\mu)=E_{\Om;A,B}(Z),\quad
Z\in\BC^{(m,n)}\end{equation*} for all $\l,\mu\in\BZ^{(m,n)}.$
Thus $E_{\Om;A,B}$ can be regarded as a function on $A_{\Om}.$
\vskip 0.1cm \end{lemma}
\begin{proof}
We write $\Om=X+iY$ with real $X,Y.$ For any
$\l,\mu\in\BZ^{(m,n)},$ we have
\begin{align*}
E_{\Om;A,B}(Z+\l\Om+\mu)&=E_{\Om;A,B}((U+\l X+\mu)+i(V+\l Y))\\
&=e^{ 2\pi i \left\{\,\s\,(\,^t\!A(U+\l
X+\mu))+\,\s\,((B-AX)Y^{-1}\,^t\!(V+\l Y))\,\right\} }\\
&=e^{ 2\pi i \left\{\,\s\,(\,^t\!AU+\,^t\!A\l X+\,^t\!A\mu)
+\,\s\,((B-AX)Y^{-1}\,^tV+B\,^t\l-AX\,^t\l) \right\} }\\
&=e^{2\pi i \left\{\,\s\,(\,^t\!AU)\,+\,\s\,((B-AX)Y^{-1}\,^tV)\right\} }\\
&=E_{\Om;A,B}(Z).\end{align*} Here we used the fact that
$^t\!A\mu$ and $B\,^t\l$ are integral. \end{proof}
\newcommand\AO{A_{\Omega}}
\newcommand\Imm{\text{Im}}
\begin{lemma}
The metric
$$ds_{\Om}^2=\s\left((\textrm{Im}\,\Om)^{-1}\,\,^t(dZ)\,d{\overline Z})\,\right)$$
is a K{\"a}hler metric on $A_{\Om}$ invariant under the action
(5.15) of $\G^J=Sp(n,\BZ)\ltimes H_{\BZ}^{(m,n)}$ on $(\Om,Z)$
with $\Om$ fixed. Its Laplacian $\Delta_{\Om}$ of $ds_{\Om}^2$ is
given by
\begin{equation*}
\Delta_{\Om}=\,\s\left( (\Imm\,\Omega)\,{ {\partial}\over
{\partial {Z}} }\,^t\!\left(  {{\partial}\over {\partial
{\overline Z}}} \right)\,
 \right). \end{equation*}
\end{lemma}

\begin{proof} The proof can be found [20].

\end{proof}

  \vskip 0.1cm We
let $L^2(\AO)$ be the space of all functions $f:\AO\lrt\BC$ such
that
$$||f||_{\Om}:=\int_{\AO}|f(Z)|^2dv_{\Om},$$
where $dv_{\Om}$ is the volume element on $\AO$ normalized so that
$\int_{\AO}dv_{\Om}=1.$ The inner product $(\,\,,\,\,)_{\Om}$ on
the Hilbert space $L^2(\AO)$ is given by
\begin{equation*}
(f,g)_{\Om}:=\int_{\AO}f(Z)\,{\overline{g(Z)} }\,dv_{\Om},\quad
f,g\in L^2(\AO).\end{equation*}
\begin{theorem}
The set $\left\{\,E_{\Om;A,B}\,|\ A,B\in\BZ^{(m,n)}\,\right\}$ is
a complete orthonormal basis for $L^2(\AO)$. Moreover we have the
following spectral decomposition of $\Delta_{\Om}$:
$$L^2(\AO)=\oplus_{A,B\in \BZ^{(m,n)}}\BC\cdot E_{\Om;A,B}.$$
\end{theorem}

\noindent $\textrm{Proof.}$ The complete proof can be found in
[20] \hfill $\square$

\vskip 1cm
\bibliography{central}



\end{document}